\documentclass[11pt]{article}
\usepackage{}
\usepackage{graphicx}
\usepackage{graphicx,subfigure}
\usepackage{epsfig}
\usepackage{amssymb,pict2e,color}
\usepackage{mathrsfs}
\usepackage{amsmath}
\usepackage{ifpdf}

\newtheorem{theorem}{Theorem}[section]

\newtheorem{lemma}[theorem]{Lemma}
\newtheorem{corollary}[theorem]{Corollary}

\newtheorem{Main results}[theorem]{Main results}
\def\whitebox{{\hbox{\hskip 1pt
 \vrule height 6pt depth 1.5pt
 \lower 1.5pt\vbox to 7.5pt{\hrule width
    3.2pt\vfill\hrule width 3.2pt}%
 \vrule height 6pt depth 1.5pt
 \hskip 1pt } }}
\def\qed{\ifhmode\allowbreak\else\nobreak\fi\hfill\quad\nobreak
     \whitebox\medbreak}
\newcommand{\proof}{\noindent{\it Proof.}\ }

\DeclareMathOperator {\gp} {gp}
\newcommand{\ignore}[1]{}

\begin {document}
\rule{0cm}{1cm}
\begin{center}
{\Large\bf On the general position set of two classes of graphs
}
\end{center}
 \vskip 2mm \centerline{Yan Yao, Mengya He, Shengjin Ji
 \footnote{ Corresponding author.\\ E-mail addresses:   yaoymath@163.com, ml19811737859@163.com, jishengjin2013@163.com,lig@sdut.edu.cn.}, Guang Li}

\begin{center}
 School of Science, Shandong University of Technology,
\\ Zibo, Shandong 255049, China
\end{center}

\begin{abstract}
The general position problem is to find the cardinality of a largest vertex subset $S$ such that no triple of vertices of $S$ lie on a common geodesic. For a connected graph $G$, the cardinality of $S$ is denoted by $\gp(G)$ and called $\gp$-number (or general position number) of $G$. In the paper, we obtain an upper bound and a lower bound regarding gp-number in all cactus with $k$ cycles and $t$ pendant edges. Furthermore, the $\gp$-number of wheel graph is determined.

\end{abstract}

\section{Introduction}

 \hspace{0.5em}

In the paper, all graphs are undirected, finite and simple. Assume that $G=(V,E)$ is a simple connected graph with vertex set $V(G)$ and edge set $E(G)$.  Let $v,r\in V(G)$. $d_G(v,r)$ denotes number of edges on a shortest $(v,r)$-path in $G$. A $(v,r)$-path with length $d_G(v,r)$ is regarded as a $v,r$-geodesic. The interval $I_G(v,r)$ of $G$ is the set of vertices $u$ so that there exists a $v,r$-geodesic which contains $u$. We refer the readers to \cite{2} for undefined terminology and notations.

The classical Dudeney's no-three-in-line problem \cite{3,4,6} is to determine the largest vertex number that can be placed in the $m\times m$ grid such that no three vertices lie on a line. The problem has been further studied in several recent papers \cite{8,10,12,13}. Later, the problem was in discrete geometry extended to General Position Subset Selection Problem \cite{5,11}, which is to obtain  a maximum subset of vertices in general position.  In \cite{5}, the  General Position Subset Selection Problem was shown to be NP-hard.

Motivated by the above two problems, the general position problem was introduced in \cite{9}. A subset $R$ of $V(G)$ is a general position set in the graph $G$ if no three vertices of $R$ lie on a common geodesic. A largest general position set is called a gp-set of $G$. The cardinality of a $\gp$-set is called the general position number ($\gp$-number for short) of $G$ and denoted by $\gp(G)$. The general position problem is to determine the $\gp$-number in the graph $G$.

In \cite{1},  Anand and Chandran et al. deduced a formula for the general position number of the complement of an arbitrary bipartite graph. In \cite{7}, Klav\v{z}ar and Yero proved that $\gp(G)\geq \omega(G_{SR})$, where $G_{SR}$ is the strong resolving graph of a connected graph $G$, and $\omega (G_{SR})$ is its clique number. In \cite{9},  Manuel and  Klav\v{z}ar determined some upper bounds on $\gp(G)$ and showed that the general position problem is NP-complete. For more properties of the general position problem, please see \cite{18k,17p}. Since researchers are very concerned about the cactus and wheel graphs, see \cite{15,14,16} for examples, it is interesting to study the bounds on the $\gp$-number for cactus and wheel graphs.

A connected graph $G$ is called a cactus if its any block is either a cycle or a cut edge, and any two cycles have no common edges. A chain cactus is a cactus graph and its each block has at most two cut vertices and each cut vertex is shared by exactly two blocks. A cycle of $G$ has one cut vertex, we call it an \emph{end-block}. A vertex of a cycle is said to be nontrivial if its degree is at least 3. Let $v$ is a cut vertex of a cycle of $G$. If a component of $G-v$ contains an end-block $C_0$(If $v$ belongs to some cycle $C_\ast$, then the component does not contain vertices of $C_\ast$), then $C_0$ is defined as an end-block of $v$. The path connects two cycles is called a \emph{cyclic path}. Let $\mathcal{C}_n^{t,k}$ be the class of all cacti of order $n$ with $k$ cycles and $t$ pendant edges. Let $\mathcal{C}_n^{k}$ be the set of all cacti of order $n$ with $k$ cycles. The wheel graph $W_n$ is the graph obtained from the cycle $C_n$ of order $n$ by adding a new vertex and connecting it to all vertices of $C_n$.

In the paper, we determine an upper bound and a lower bound of $\gp$-number for cacti graphs, and also characterize the property of the extremal graphs that attain the bounds. In addition, the $\gp$-number is obtained for wheel graphs.

\section{The upper bound of cacti regarding $\gp$-number}

 From \cite{9}, we know that $\gp(C_3)=3$, $\gp(C_4)=2$ and $\gp(C_n)=3$ for $n\geq 5$.
\begin{lemma}\label{bc1}
 Let graph $G$ be a cactus with $k(\geq 2)$ cycles, and $S$ be a $\gp$-set of $G$. If $C_0$ is a cycle of $G$, then $|V(C_0) \cap S |<3$.
\end{lemma}\par

\begin{proof}
Let $C_1$ and $C_2$ be two cycles of $G$ connected by a path $P$ (maybe trivial). Assume that $V(C_1)\cap V(P)=u$ and $V(C_2)\cap V(P)=r$. Let $u_1$ and  $u_2$, $u_3$ and  $u_4$ be the adjacent vertices of $u$ and $r$ in $C_1$ and $C_2$, respectively. Note that $R=\{u_1,u_2,u_3,u_4\}$ is a general position set of $G$. Therefore, $\gp(G)\geq 4$.

Suppose $|V(C_0) \cap S |= 3$ and let $V(C_0) \cap S=\{x, y, z\}$. This means that $I(x,y)\cup I(y,z) \cup I(z,x)=V(C_0)$ and $I(x,y)\cap I(y,z)=y $, $I(y,z)\cap I(z,x)=z $, $I(x,y)\cap I(z,x)=x $ . Since $G$ is a cactus, the cycle $C_0$ possesses at least one nontrivial cut vertex, marked as $v$. Let $H_0$ be a nontrivial subgraph of $G$ such that $V(H_0) \cap V(C_0)=v$. Without loss of generality, assume that $v \in I(x,y) $. For some $v_0 \in V(H_0)$, every path connecting the vertex $v_0$ and all vertices in $C_0$ goes through $v$. Therefore $x $(or  $y$) $\in I(v_0, z)$ and we deduce that $|S\cap V(H_0)|=0$. Then $\gp(G)=3$, a contradiction.
\qed
\end{proof}

\cite{9} determined the $\gp$-number of trees. We present it here.
\begin{lemma}\label{bc2}
If $L$ is the set of leaves of a tree T, then $\gp(T)=|L|$.
\end{lemma}\par
For a cactus $G$, if a tree $T_k$ and a cycle $C_t$ (or a cyclic path $P_s$) of a cactus $G$ share a common vertex $v$, we say that $T_k$ is a \emph{pendant tree} in $G$  associated  with $v$ and $v$ is the\emph{ root} of $T_k$. A root vertex $v$ is referred as a \emph{nontrivial cut vertex} of $G.$
\begin{lemma}\label{bc3}
 If a $G\in \mathcal{C}_n^{t,k}$ and $\mathcal {T}$ is the set of all pendant trees except for their roots in $G$, then there exists a $\gp$-set $S$ such that $|S\cap V(\mathcal {T})| = t $.
\end{lemma}\par

\begin{proof}
  Let $G\in \mathcal{C}_n^{t,k}$ has $w$ pendant trees labeled as  $T_1,\ldots,T_{\ell},T_{\ell+1},\ldots,T_{w}$ with roots $v_1,\linebreak\ldots,v_{\ell},v_{\ell+1},\ldots,v_{w}$, where the first $\ell$ roots (resp. the last $w-\ell$ roots) belong to cycles(resp. cyclic paths). Let $\mathcal {T}_1=\bigcup _{i=1}^{\ell}\{T_i - v_i\}$, $\mathcal {T}_2=\bigcup _{i=\ell+1}^{w}\{T_i - v_i\}$ and $\mathcal {T}=\mathcal {T}_1\cup \mathcal {T}_2$. Let $L_i$ $(1\leq i \leq w)$ be the set of leaves of the pendant tree $T_i$  and $X_i$ $(1\leq i \leq w)$ be the set of vertices with degree $1$ in $G$. Obviously, $|L_i|=|X_i|$. Let $S'$ be a $\gp$-set of $G$.

 {Claim 1.} $S'$  doesn't contain any root $v_i(\ell+1\leq i\leq w).$

  Assume that there is some $v_{i_0}\in S'$ for $\ell+1\leq i_{0}\leq w$.
As shown in Fig.2, $v_{i_0}$ is a nontrivial cut vertex of $G$, which, besides $T_{i_0}$, has $p(\geq 2)$ components containing cycles labeled as $H_1,H_2,\ldots,H_{p}$. For every triple including $v_{i_0}$ of $S'$ marked as $\{x,y,v_{i_0}\}$, we conclude that $x$ and $y$ simultaneously belong to some component (or the pendant tree) of $v_{i_0}$, assume that $H_1$. If not, suppose that $x\in V(H_j)$ and $y\in V(H_k)$(or $V(T_{i_0}-v_{i_0}$)), then  any $x,y$-geodesic is via $v_{i_0}$, contradicted with $S'.$  By the way, we get a new general position set $S={S'}\cap V(H_1)\cup \{u_2,u_3,\cdots,u_{p},u_{p+1}\}.$ Obviously, $|S|\geq |S'|+2$, which is a contradiction.

We now consider the relationship between $S'$  and the remaining roots of $G$.

 If $S'$ doesn't contain any root $v_{i}(1\leq i\leq \ell)$, we take $S=S'$. Every vertex with degree $1$ can not lies on the geodesics of any other pair pendant vertices. Hence, $|S\cap V(\mathcal {T}_1)|=|\sum^{\ell}_{i=1}S\cap (V(T_i)-v_i)|=\sum^{\ell}_{i=1}|S\cap V(T_i)|=\sum^{\ell}_{i=1}|X_i|=\sum^w_{i=1}|L_i|$ by Lemma \ref{bc2}.

 We hence suppose that $S'$ includes some root $v_{i_0}$ which belongs to a cycle $C_0$. Evidently, $S' = ( S'\cap V(G-T_{i_0})) \cup ( S'\cap (V(T_{i_0})\setminus \{v_{i_0}\})) \cup \{v_{i_0}\}$. We claim that $T_{i_0}$ is a path. Note that every path from the vertices in $V(G-T_{i_0})$ to the vertices of $V(T_{i_0})\setminus \{v_{i_0}\}$ is via $v_{i_0}$, which implies $S'\cap V(G-T_{i_0})=\emptyset$ or $S'\cap (V(T_{i_0})\setminus \{v_{i_0}\})=\emptyset$.  Firstly, suppose that  $S'\cap (V(G-T_{i_0}) ) =\emptyset$. Mark the two adjacent vertices of $v_{i_0}$ as $x$ and $y$ in $C_0$. Then we obtain a new general position set $S''=\{x, y\} \cup (S'\cap (V(T_{i_0})\setminus \{v_{i_0}\}))$. this means that  $|S''|> |S'|$, which is
a contradiction with the maximum of $S'$. That is to say, $S'\cap (V(T_{i_0})\setminus \{v_{i_0}\}))= \emptyset$. If $|X_{i_0}|\geq 2$, we have a new general position set $S'''=(S'\cap V(G-T_{i_0}) )\cup X_{i_0}$. Then $|S'''|> |S'|$,  contradicted with the maximum of $S'$. So $|X_{i_0}|=1$ which implies that $T_{i_0}$ is a path. Let $S=(S' \setminus \{v_{i_0}\})) \cup X_{i_0}$. Clearly, $|S|=|S'|$ and $S$ is a $\gp$-set of $G$.

  Assume now that $S'$ contains some of these roots. Then, repeating the above process, we obtain a new $\gp$-set $S$ of $G$ such that all roots are not its elements.

By Lemma \ref{bc2}, we conclude that $|S\cap V(\mathcal {T})|= \sum\limits_{i=1}^{w}|S \cap V(T_i)|=\sum\limits_{i=1}^{w}|X_i|=\sum\limits_{i=1}^{w}|L_i|=t$.
\qed
\end{proof}

 Let $\mathcal {C}$ be the set of all cycles and $\mathcal {T}$ is the set of all pendant trees except for their roots of a cactus $G$. From Lemmas \ref{bc1}, \ref{bc3}, we deduce $\gp(G)= |S\cap V(\mathcal {C})|+|S\cap V(\mathcal {T})|\leq 2k+t$. In order to check that the bound is best, we next show the structure of all cacti such that their $gp$-number equal $2k+t$.

Let $G\in \mathcal{C}_n^{t,k}$ and $C_{\ell}=u_1u_2 \ldots u_{\ell}u_1$ be a cycle of $G$. For two vertices $u_i$ and $u_j$, clearly, $C_{\ell}$ is consist of two $(u_i,u_j)$-paths. If $u_i$ and $u_j$ are two vertices and all others cut-vertices of $C_{\ell}$ belong to one $(u_i,u_j)$-path. Then this path is referred to as a
\emph{$(u_i,u_j)$-cut-path} and use $d_c(u_i,u_j)$ to denote the number of edges in the $(u_i,u_j)$-cut-path. Let $D_c$ be $\min_{u_i,u_j\in V(C_{\ell})}d_c(u_i,u_j)$ of $C_{\ell}.$ We call cycle $C_{\ell}$ a good cycle if $D_c$ is no more than $\frac{\ell}{2}$ for even $\ell$ or $\lfloor\frac{\ell}{2}\rfloor-1$ for odd $\ell$. Trivially, the cycle as an end-block of cactus is a good cycle. Moreover, if $D_c$ is at least $\frac{\ell}{2}-1$ for even $\ell$ and $\lfloor\frac{\ell}{2}\rfloor$ for odd $\ell$, then the cycle is regarded as a bad cycle. Especially, if $C_{\ell}$
has exactly two cut-vertices, then
$C_{\ell}$ is a bad cycle, if $D_c=\frac{\ell}{2}$ for even $\ell$, otherwise, $D_c=\lfloor\frac{\ell}{2}\rfloor.$
Let $H_0=v_1v_2 \ldots v_{\ell}v_1$ be a cycle of $G_0$ (see Fig.1). Assume that $H_i(1\leq i\leq \ell)$ is a connected subgraph with $V(H_i)\cap V(H_0)=v_i$ and $|V(H_i)|=n_i$ for $i\in \{0,1,\ldots ,\ell\}$.

\vspace{-7mm}
\begin{center}
\begin{picture}(288.2,119.1)\linethickness{0.8pt}

\cbezier(17.2,80.8)(18.7,69.6)(63.1,69.6)(64.7,80.8)
\cbezier(17.2,80.8)(18.7,92.1)(63.1,92.1)(64.7,80.8)
\put(50,78){\footnotesize$T_{i_0}$}
\put(45.5,70){\footnotesize$v_{i_0}$}
\put(8.7,48.9){\tiny$H_1$}
\put(27.8,47.8){\tiny$H_2$}
\put(57.8,45.8){\tiny$H_{p}$}
\put(75,45){$G$}
\put(34.7,51.8){\circle{15.6}}
\put(64.2,49.8){\circle{16.6}}
\put(15.7,52.8){\circle{17}}
\Line(39.5,72.9)(20.2,60.3)
\Line(39.5,72.9)(38.2,58.8)
\Line(39.5,72.9)(61.7,57.8)
\Line(31.2,83.3)(39.5,72.9)
\put(60,62){\footnotesize$u_{p}$}
\put(38,62){\footnotesize$u_2$}
\put(14,63){\footnotesize$u_1$}
\put(20,78){\footnotesize$u_{p+1}$}
\put(31.2,83.3){\circle*{4}}
\put(39.5,72.9){\circle*{4}}
\put(20.2,60.3){\circle*{4}}
\put(38.2,58.8){\circle*{4}}
\put(61.7,57.8){\circle*{4}}
\put(50.2,56){\circle*{2}\color{black}\circle*{1}}
\put(55.2,56){\circle*{2}\color{black}\circle*{1}}
\put(45.7,56){\circle*{2}\color{black}\circle*{1}}


\Line(195.2,91.9)(226.2,91.9)
\Line(195.2,91.9)(178.7,69.9)
\Line(178.7,69.9)(196.2,45.9)
\Line(196.2,45.9)(224.2,45.9)
\put(224.2,45.9){\circle*{4}}

\Line(226.2,91.9)(237.7,71.4)
\put(226.2,91.9){\circle{0}}
\put(246,71){\circle{16}}
\put(187,47){\circle{16}}
\put(171,69){\circle{16}}
\put(188.2,95.4){\circle{16}}
\put(234,95.9){\circle{16}}
\put(225,95){\tiny
$H_1$}
\put(242,67.5){\tiny
$H_2$}
\put(177,44){\tiny
$H_{\ell-2}$}
\put(157,64.5){\tiny
$H_{l-1}$}
\put(184,95){\tiny
$H_{\ell}$}
\put(202,66.5){\footnotesize
$H_0$}
\put(219.2,81.4){\footnotesize
$v_1$}
\put(224.7,68.4){\footnotesize
$v_2$}
\put(194,49.9){\footnotesize
$v_{\ell-2}$}
\put(182.7,65.9){\footnotesize
$v_{\ell-1}$}
\put(195.7,81.9){\footnotesize
$v_{\ell}$}
\put(195.2,91.9){\circle*{4}}
\put(226.2,91.9){\circle*{4}}
\put(178.7,69.9){\circle*{4}}
\put(196.2,45.9){\circle*{4}}

\put(234.7,61.9){\circle*{2}}
\put(231,57){\circle*{2}}
\put(228,52){\circle*{2}}

\put(226.2,91.9){\circle*{4}}
\put(237.7,70.9){\circle*{4}}

\put(178.2,69.9){\circle*{4}}
\put(195.2,91.9){\circle*{4}}
\put(226.2,92.4){\circle*{4}}

\put(255,40){$G_0$}
\end{picture}
\end{center}\vspace{-15mm}
\centerline{Fig.1 The two graphs $G$, $G_0$ used in Lemma \ref{bc3} and Theorem \ref{bc4}.}\vspace{3mm}

\begin{theorem}\label{bc4}
Let $G_0 \in\mathcal{C}_n^{t,k}$. If all cycles of $G_0$ are good, then $\gp(G_0)=2k+t$.
\end{theorem}\par

\begin{proof}
 Let $G_0 \in\mathcal{C}_n^{t,k}$ and $H_0$ is a cycle of $G_0$. Suppose that $S$ is a $\gp$-set of $G_0$. In order to verify $|S|=2k+t$, it is sufficient to show $H_0$ is good. If $H_0$ is an end-block of $G_0$, then it is easy to imply that $|V(H_0)\cup S|=2$ and $H_0$ is a good cycle. We now assume that $H_0$ is not an end-block as shown in Fig.1.(It includes at least two cut vertices.)

 Assume that $H_0$ is not a good cycle. Hence its $D_c$ is no less than $\frac{\ell}{2}-1$ for even $\ell$ or $\lfloor\frac{\ell}{2}\rfloor$ for odd $\ell$. Let $v_{i_0}$ and $v_{j_0}$ be two cut-vertices of $H_0$ for which $D_c=d_c(v_{i_0},v_{j_0})$. We first consider a special case that
$H_0$ contains only two cut vertices. Then, $D_c\leq\lfloor\frac{\ell}{2}\rfloor$, and we deduce that $|S\cap V(H_0)|\leq 1.$

 Suppose that there are at least three cut vertices in $V(H_0)$. If $d_c(v_{i_0},v_{j_0})\geq \lfloor\frac{\ell}{2}\rfloor+1$, then there exists a cut vertex $v_r$ which belongs to the $(v_{i_0},v_{j_0})$-cut-path of the cycle $H_0$, such that one vertex is not on the geodesic of the other two vertices in the triple $\{v_{i_0},v_{j_0},v_r\}.$ Let now $v_{k_0}(\neq v_{i_0},v_{j_0},v_r)$ be an arbitrary vertex of $H_0$. Assume that $v_{k_0}$ is one vertex of the shortest $(v_{i_0},v_{j_0})$-path. Then, It is easy to deduce that every shortest path form the vertices in $H_{v_{i_0}}$ to the vertices in $H_{v_{j_0}}$ goes through $v_{k_0}$. Consequently, at least one of $|S \cap V(H_{v_{i_0}})|$, $|S \cap V(H_{v_{j_0}})|$ and $|S \cap V(H_{v_{r}})|$ is equal to zero, and $\gp(G_0)<2k+t$ by Lemma \ref{bc1} and \ref{bc3}.

 Assume that $d_c(v_{i_0},v_{j_0})\geq \lfloor\frac{\ell}{2}\rfloor$. If $\ell$ is even, then $ \lfloor\frac{\ell}{2}\rfloor=\frac{\ell}{2}$.
For a vertex $u$ of $H_{v_{i_0}}-v_{i_0}$ and a vertex $w$ of $H_{v_{j_0}}-v_{j_0}$, all vertices of $H_0$ lie on some $(u,w)$-geodesic. Hence $|S\cap V(H_0)|=0$ or $|S\cap V(H_{v_{i_0}})|=0$ or $|S\cap V(H_{v_{j_0}})|=0$. According to Lemma \ref{bc1} and Lemma \ref{bc3}, $\gp(G_0)<2k+t.$ If $\ell$ is odd, then $ \lfloor\frac{\ell}{2}\rfloor=\frac{\ell-1}{2}$. Using a similar way for even $\ell$, it is verified that $\gp(G_0)<2k+t.$

Suppose $d_c(v_{i_0},v_{j_0})= \lfloor\frac{\ell}{2}\rfloor-1=
\frac{\ell}{2}-1$ for even $\ell$. We conclude that any vertex which lies on a $(v_{i_0},v_{j_0})$-cut-path is not a element of $S$. Let $P_0$ be another $(v_{i_0},v_{j_0})$-path in $H_0$. Clearly, $|S\cap(V(P_0)\backslash\{v_{i_0},v_{j_0}\})|=1$,
which implies $\gp(G)<2k+t$.

Therefore, we next consider that $H_0$ is a good cycle.

Suppose $d_c(v_{i_0},v_{j_0})\leq \frac{\ell}{2}-2$ for even $\ell$. We take the two vertices which are adjacent to $v_{i_0}$ and $v_{j_0}$ on another $(v_{i_0},v_{j_0})$-path in $H_0$, and mark them as $(v'_{i_0},v'_{j_0})$, respectively. It is clear that the two vertices have degree two. For an arbitrary vertex $u_0$ in $G-V(H_0)$, $v'_{i_0}$ is not on a $(u_0,v'_{j_0})$-shortest path and
$v'_{j_0}$ is also not on a $(u_0,v'_{i_0})$-shortest path. If $|S\cap V(H_0)|<2$, then let $S'$ be the set obtained from $S$ by removing all the elements of its subset $S\cap V(H_0)$ and adding the two vertices $v'_{i_0}$ and $v'_{j_0}$. We deduce that $S'$ is a general position set with $|S'|\geq |S|+1$, which contradicts the maximum of $S$. Hence, $|S\cap V(H_0)|=2$.

Suppose $d_c(v_{i_0},v_{j_0})\leq \lfloor\frac{\ell}{2}\rfloor-1=
\frac{\ell-1}{2}-1$ for odd $\ell$. Let $v'_{i_0}$ and $v'_{j_0})$ denote the two vertices of degree two which are adjacent to $v_{i_0}$ and $v_{j_0}$ on another $(v_{i_0},v_{j_0})$-path, respectively. Hence, for an arbitrary vertex $w_0$ in $G-V(H_0)$, $G-V(H_0)$, $v'_{i_0}$ is not on a $(w_0,v'_{j_0})$-shortest path and
$v'_{j_0}$ is also not on a $(w_0,v'_{i_0})$-shortest path. By a similar way as above, we get
$|S\cap V(H_0)|=2$.

Therefore, we obtain that if $H_0$ is a good cycle, there is a $\gp$-set $S$ for which $|S\cap V(H_0)|=2$.

Together with Lemma \ref{bc3}, the proof is complete.
\qed
\end{proof}

\begin{center}
\begin{picture}(325.7,83.8)\linethickness{0.8pt}

\Line(72.7,73.9)(55.2,73.9)
\Line(72.7,73.9)(40.7,48.9)
\Line(55.2,73.9)(40.7,48.9)
\Line(40.7,48.9)(21.2,25.4)
\Line(40.7,48.9)(59.2,24.4)
\Line(40.7,48.9)(31.7,24.4)
\put(43.8,68){\footnotesize$k$}
\put(13.9,15){\footnotesize $n-2k-1$}

\put(70,25){$B_0$}
\Line(15.2,73.4)(32.7,73.9)
\Line(15.2,73.4)(40.7,48.4)
\Line(32.7,73.9)(41.2,48.9)

\put(40.7,48.9){\circle*{4}}
\put(72.7,73.9){\circle*{4}}
\put(55.2,73.9){\circle*{4}}
\put(21.2,25.4){\circle*{4}}
\put(31.7,24.4){\circle*{4}}

\put(40,24.4){\circle*{2}}
\put(45,24.4){\circle*{2}}
\put(50,24.4){\circle*{2}}

\put(59.2,24.4){\circle*{4}}

\put(15.2,73.4){\circle*{4}}
\put(32.7,73.9){\circle*{4}}
\put(40.7,48.4){\circle*{4}}
\put(41.2,48.9){\circle*{4}}
\put(38,74.4){\circle*{2}}
\put(43,74.4){\circle*{2}}
\put(48,74.4){\circle*{2}}

\put(150.2,57.7){\circle{28.5}}
\put(142.5,58.4){\footnotesize$C_1$}
\put(184.7,59.8){\circle{27.7}}
\put(169.3,35.9){\circle{28.2}}
\Line(140,47.8)(125.7,46.8)
\Line(155.2,36.2)(143.2,35.3)
\Line(143.2,35.3)(129.7,36.3)
\Line(143.2,35.3)(131.2,27.3)
\put(180,56.6){\footnotesize$C_2$}
\put(259.7,41.3){\circle{20.4}}
\put(276.2,58.8){\circle{22.9}}
\put(306.2,52.2){\circle{27.5}}
\Line(315.3,41.9)(314.2,28.3)
\Line(315.3,41.9)(325.7,31.8)

\put(190,25){$B_1$}
\put(330,25){$B_2$}
\put(277.9,29.3){\footnotesize$C_1$}
\put(251.8,39.5){\footnotesize$C_2$}
\put(163.4,30.4){\footnotesize$C_3$}

\put(272.4,56.9){\footnotesize$C_3$}
\put(299.9,51.4){\footnotesize$C_4$}
\put(284.2,33.3){\circle{30.6}}
\Line(297.3,25.4)(297.3,25.4)

\put(151.9,71.9){\circle*{4}\color[rgb]{1,0.98,0.94}\circle*{3}}
\put(136.6,62){\circle*{4}\color{black}\circle*{3}}
\put(198.5,58.6){\circle*{4}\color{black}\circle*{3}}
\put(171.9,65.2){\circle*{4}\color{black}\circle*{3}}
\put(159.7,46.2){\circle*{4}\color[rgb]{1,0.98,0.94}\circle*{3}}
\put(164.1,61.2){\circle*{4}\color{black}\circle*{3}}
\put(160.5,24.9){\circle*{4}\color{black}\circle*{3}}
\put(140,47.8){\circle*{4}\color[rgb]{1,0.98,0.94}\circle*{3}}
\put(125.7,46.8){\circle*{4}\color{black}\circle*{3}}
\put(155.2,36.2){\circle*{4}\color[rgb]{1,0.98,0.94}\circle*{3}}
\put(143.2,35.3){\circle*{4}\color[rgb]{1,0.98,0.94}\circle*{3}}
\put(129.7,36.3){\circle*{4}\color{black}\circle*{3}}
\put(131.2,27.3){\circle*{4}\color{black}\circle*{3}}
\put(172.6,22.2){\circle*{4}\color[rgb]{1,0.98,0.94}\circle*{3}}
\put(181.2,28.3){\circle*{4}\color[rgb]{1,0.98,0.94}\circle*{3}}
\put(177.5,47.9){\circle*{4}\color[rgb]{1,0.98,0.94}\circle*{3}}
\put(183,39.1){\circle*{4}\color{black}\circle*{3}}

\put(256.7,51.3){\circle*{4}\color{black}\circle*{3}}
\put(255.2,33.3){\circle*{4}\color{black}\circle*{3}}
\put(265.3,55.2){\circle*{4}\color{black}\circle*{3}}
\put(287.3,61.7){\circle*{4}\color{black}\circle*{3}}
\put(272.6,69.6){\circle*{4}\color[rgb]{1,0.98,0.94}\circle*{3}}
\put(294,58.4){\circle*{4}\color{black}\circle*{3}}
\put(307.1,65.9){\circle*{4}\color[rgb]{1,0.98,0.94}\circle*{3}}
\put(319.7,55.1){\circle*{4}\color{black}\circle*{3}}
\put(315.3,41.9){\circle*{4}\color[rgb]{1,0.98,0.94}\circle*{3}}
\put(314.2,28.3){\circle*{4}\color{black}\circle*{3}}
\put(325.7,31.8){\circle*{4}\color{black}\circle*{3}}

\put(279.1,47.8){\circle*{4}\color[rgb]{1,0.98,0.94}\circle*{3}}
\put(295,44.2){\circle*{4}\color[rgb]{1,0.98,0.94}\circle*{3}}
\put(269.7,39.8){\circle*{4}\color[rgb]{1,0.98,0.94}\circle*{3}}
\put(270.7,27.8){\circle*{4}\color{black}\circle*{3}}
\put(277.7,19.5){\circle*{4}\color[rgb]{1,0.98,0.94}\circle*{3}}
\put(289.3,18.9){\circle*{4}\color[rgb]{1,0.98,0.94}\circle*{3}}
\put(297.3,25.4){\circle*{4}}
\put(297.3,25.4){\circle*{4}\color[rgb]{1,0.98,0.94}\circle*{3}}
\put(299.4,34.8){\circle*{4}\color{black}\circle*{3}}

\put(269.7,39.3){\circle*{4}\color[rgb]{1,0.98,0.94}\circle*{3}}
\end{picture}
\end{center}\vspace{-5mm}
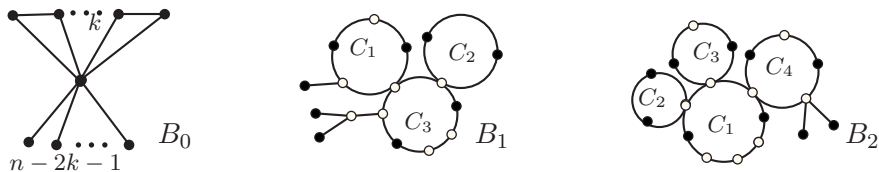
\centerline{Fig.2 The three graphs used in Theorem \ref{bc5} and for two examples}
\vspace{3mm}

    By means of Theorem \ref{bc4} we deduce that there are a lot of graphs belonging to $\mathcal{C}_n^{t,k}$ for which their $\gp$-number equal to $2k+t$ and these graphs have different shapes.  

    Let $A_1$ be a cactus in which all cycles and all pendent trees share a common cut vertex $v$. Let $R$ be the neighbors of $v$ in cycles and the pendant vertices in $A_1$. In fact, $R$ is a general position ,and $|R|=2k+5$. Furthermore, $\gp(A_1)\geq 2k+t$ by Lemma \ref{bc1} and Lemma \ref{bc3}, which implies $\gp(A_1)=2k+t$.

In addition, let $A_2$ denote a chain cactus with the following two properties: (a) the length of every cycle in $A_2$ is at least $5$ and its two cut vertices (if they exist) are adjacent (b) there is one cut vertex of some cycle possessing pendant tree. Evidently, all cut vertices of cycles and all pendant vertices form a general position set of $A_2$, denote it by $R'$. Then $|R'|=2k+t$ and $\gp(A_2)\leq 2k+t.$ So $\gp(A_2)=2k+t$.

We now give two concrete examples.
 Let $B_1\in \mathcal{C}_n^{t,k}$ with $k=3$, $t=3$, $n=17$, and $B_2\in \mathcal{C}_n^{t,k}$ with $k=4$, $t=2$, $n=19$, see Fig.2. We have their general position sets $S$ and $S'$ that consists of all solid vertices, respectively. Evidently, $|S|=9=2k+t$ and $|S'|=10=2k+t$.

By means of Theorem \ref{bc4}, the main result holds.

\begin{theorem}\label{bc5}
 If $G\in \mathcal{C}_n^{t,k}$, then
 \begin{equation*}
   \gp(G)\leq \max\{3, 2k+t\}.
 \end{equation*}
 with the equality holds if and only if all cycles of $G$ are good.
\end{theorem}

\begin{proof}
{\bf Case 1.} $k=1$ with at most one pendant edge.

It is easy to verify that $\gp(G)=3$.

{\bf Case 2.} $k\geq2$ or $k=1$  with at least two pendant edges.

We conclude that $\gp(G)\leq 2k+t$ by Lemma \ref{bc1}, Lemma \ref{bc3}. Moreover, Theorem \ref{bc4} implies that $\gp(G)= 2k+t$ if and only if every cycle of $G$ is good.
\qed
\end{proof}

\begin{corollary}\label{bc6}
For any $G\in \mathcal{C}_n^{k}$, we have
 \begin{equation*}
   \gp(G)\leq \max \{3,n-1\}.
 \end{equation*}
 with equality if and only if $G\cong C_3$ or $G\cong B$.
\end{corollary}
\begin{proof}
{\bf Case 1.} $G$ is a $C_3$.

$gp(G)=3$.

  {\bf Case 2.} $G$ is not a $C_3$.

Let $t$ be the number of pendent edges in $G$. It is easy to verify that   $t\leq n-2k-1$. Because the $k$ is limited , we can get $\gp(G)\leq 2k+t \leq 2k+n-2k-1=n-1$.

 The quality $\gp(G)=n-1$ means that $t=n-2k-1$. That is to say, every cycle of $G$ triangle  and they share a common vertex, we have $G\cong B_0$ (see Fig.2). Actually, it is easy to check that all vertices except for the cut vertex form a general position set of $G$. Hence, $\gp(B_0)=n-1.$ \qed
\end{proof}
\section{A lower bound for cacti}
In the section, we will prove some lower bounds for cacti. Note that for a graph $G\in\mathcal{C}_n^{t,k}$, its pendant vertices consist of a general position set. Hence there exist an obvious lower bound, $\gp (G)\geq t$. It is interesting to determine all graphs attaining this bound. In addition, we also consider cactus graph with $t=0.$

If $G$ does not contain pendant vertices ($t=0$), we have the following result.

\begin{theorem}\label{lbound}
If $G\in \mathcal{C}_n^{0,k}$ has $k_1(\geq 1)$ odd cycles and $k-k_1$ even cycles, then $\gp(G)\geq k_1+2$ for $k_1\geq 2$, otherwise, $\gp(G)=4$. Equality holds if and only if $G$ is a chain cactus for which, except for two end-blocks, $D_c$ of every even cycle equals to the half number of its vertices and $D_c$ of every odd cycle equals to the floor of the half number of its vertices.
\end{theorem}

\proof Let $G$ be a cactus from $\mathcal{C}_n^{0,k}$ having the minimal $\gp$-number, and let $S$ be a $\gp$-set of $G$. We next claim that $G$ has just two end-blocks. Otherwise, suppose that $G$ possesses at least three end-blocks. Two cases should be discussed.

{Case 1} There is a vertex $v$ which does not lie any cycle of $G$ has at least three end-blocks.

Mark the neighbors of $v$ as $v_1,v_2,\ldots,v_{t_1}(t_1\geq 3)$. Assume that $v_1,v_2$ and $v_3$ have the three end-blocks $C_1,C_2$ and $C_3$ and their cuts writes as $u_1,u_2$ and $u_3$, respectively. Let $u$ be a vertex of $C_3$ such that $d(u_3,u)=\frac{|C_3|}{2}$ for even $C_3$ and $d(u_3,u)=\lfloor\frac{|C_3|}{2}\rfloor$ for odd $C_3$. Let $G'$ be the graph from $G$ by deleting the edge $vv_1$ and joining $u$ and $v_1$. Set $S'=S\setminus(S\cap C_3)(resp.(S\setminus(S\cap C_3))\cup\{w\})$ for even (resp. odd) $C_3$, where $w$ is a internal vertex on the longer $(u_3,u)$-path of $C_3$. Obviously, $S'$ is a $\gp$-set of $G'$ and $|S\cap V(C_3)|=2$. Hence, $|S'|=|S|-2(\,or\,|S|-1)$ which results in a contradiction.

{Case 2} Let $v'$ be a vertex that lies on some cycle of $G$, and suppose it has at least three end-blocks.

Label $N(v')$ as $\{v'_1,v'_2,\ldots,v'_{t_2}\}(t_2\geq 3)$. Let $v'_1,v'_2$ and $v'_3$ have three end-blocks $C'_1,C'_2$ and $C'_3$. Moreover, write their unique cuts as $u'_1,u'_2$ and $u'_3$, respectively. Let $u'$ denote the vertex of $C'_3$ such that $d(u'_3,u')=\frac{|C'_3|}{2}$ for even $C'_3$ and $d(u'_3,u')=\lfloor\frac{|C'_3|}{2}\rfloor$ for odd $C'_3$. Let $G''$ be the graph obtained from $G$ by removing the edge $vv'_1$ and connecting $u'$ and $v'_1$. Let $S''=S\setminus(S\cap C'_3)(resp.(S\setminus(S\cap C'_3))\cup\{w'\})$ for even (resp. odd) $C'_3$, where $w'$ lies on the longer $(u'_3u')$-path of $C'_3$. Evidently, $S''$ is a $\gp$-set of $G'$ and $|S\cap V(C'_3)|=2$. Consequently, $|S''|=|S|-2(\,or\,|S|-1)$ which contradicts the minimum of $G$.

Note that, except for two end-blocks, all others internal cycles contain two cut vertices. Furthermore, we discuss the contribution of internal cycles to $S$, as an bad even cycle $C_0$, $|S\cap V(C_0)|=0$, and as a bad odd cycle $C_0$, $|S\cap V(C_0)|=1$. In order to keep $G$ with minimum, even cycles should be internal cycles as more as possible. As we known, $G$ has $k_1$ odd cycles and $k-k_1$ even cycles, we hence arrive at $\gp(G)=k_1-2+4$ for $k_1\geq 2$ and $\gp(G)=4$ for $k_1=1.$ \qed

Based on the above conclusion, we deduce the following two results.

\begin{corollary}\label{lboundo}
If $G\in\mathcal{C}_n^{0,k}$ has $k$ odd cycles, then $\gp(G)\geq k+2$ with equality if and only if $G$ is a chain cactus for which, except for two end-blocks, $D_c$ of every cycle equals to the floor of the half number of its vertices.
\end{corollary}

\begin{corollary}\label{lbounde}
If $G\in\mathcal{C}_n^{0,k}$ has $k$ even cycles, then $\gp(G)\geq 4$ with equality if and only if $G$ is a chain cactus for which, except for two end-blocks, $D_c$ of every cycle equals to the half number of its vertices.
\end{corollary}
\begin{theorem}\label{lboundt}
If $G\in\mathcal{C}_n^{t,k}$, then $\gp(G)\geq t$ with equality if and only if the cycles with at least three cut vertices of $G$ are bad and the cycles with two cut vertices of $G$ are bad and even.
\end{theorem}

\proof $G$ denotes an element of $\mathcal{C}_n^{t,k}$. By the property of the set of pendant vertices, $\gp(G)\geq t.$ If equality holds, then $S$ is a $\gp$-set. This means that the vertices of all cycles do not contribute to $S$. By means of the notations of good or bad cycles, equality holds if and only if $G$ contains $k$ bad cycles without these bad odd cycles having two cut vertices. \qed

\section{Wheel graphs}
In the section, we will deduce the 	gp-number of wheel graph. For convenience, we now introduce some notation. For a given graph $G$, $e_1=u_1u_2,e_2=v_1v_2\in E(G)$. The \emph{distance} of  $e_1$ and $e_2$, denoted by $d( e_1,e_2)$, is defined as $\min\{d(u_i,v_j),i,j\in\{1,2\}$. The vertex with degree $n$ in $W_n$ is called  the center of $W_n$ and denoted by $w$. Let $C_n$ be the unique cycle of $W_n-w,$  and $S$ denote a $\gp$-set of $W_n$.
  Let $G_0=C_n-S$ for short.

\begin{theorem}\label{bc6}
If $n\geq 3$, then $\gp(W_n)= \left\{\begin{array}{ll} 4, & \text{if }  n=3,\\
3, & \text{if }  n=4,5,\\
     \lfloor \frac{2}{3} n\rfloor, & \text{if }  n\geq 6 .
    \end{array}\right.$
\end{theorem}
\begin{proof}
    If $n=3$, we have $\gp(W_n)=4$. If $n=4,5$, then $\gp(W_n)=3$. We next assume that $n\geq 6$. Let $S$ be a gp-set of $W_n$. Write $w$ as the center of $W_n.$

{\bf Claim 1.} If $n\geq 6$, then the center $w \notin S$.

{\bf Proof of Claim 1.} If $w\in S$, then $|S|=3$. Since $n\geq 6$, there are two edges for which their distance no less than 2 in $C_n$, denote their ends as $x_1,x_2$ and $x_3,x_4$. Then $S'=\{x_1,x_2,x_3,x_4\}$ is a general position set, which implies that $|S'| > |S|$. So $w \notin S$.
\qed
 Note that $C_n[S]$ contains no path $P_3$. Hence, it consists of some $K_2$ and $K_1$.

{\bf Claim 2.} $G_0$ contains at most one $K_2$.

{\bf Proof of Claim 2.}  Suppose $G_0$ has at least two $K_2$. Let $e_1=v_1v_2$ be an edge in $G_0$. It is not difficult to find another edge $e_2=v_{x-1}v_x \in E(G_0)$ such that it is closest to $e_1$ in $C_n$.  We now construct a new subset $S'$ from $S$ by the following process. If the vertex  $v_i\in S$ $(3 \leq i \leq  x-2)$, then $v_{i-1}\in S'$, for others vertices of $S$, we copy them to $S'$. Since $v_{x-2} \in S$, we have $v_{x-2} \notin S'$. We add $v_{x-1}$ to $S'$. In fact, every path in $W_n$ possessing some triple of $S'$ is not geodesic.  Hence, $S'$ is a general position set of $W_n$.
Then $|S'|=|S|+1>|S|$, a contradiction.
\qed
 Let $C_n=v_1v_2 \cdots v_nv_1 $.
 Based on the above two claims, the $\gp$-set $S$ is a subset of $V(C_n)$ and $C_n[S]$ consists of some $K_1$ and $K_2$ for which every two consecutive components are  separated by $K_1$ or $K_2$(at most one). For convenience's sake, the number of $K_1$ and $K_2$ in $C_n[S]$ are $\ell_1$ and $\ell_2$, respectively. Hence, $\gp(W_n)=|S|=\ell_1+2\ell_2$. In addition, if $G_0$ contains one $K_2$, $|G_0|=\ell_1+\ell_2+1$; otherwise, $|G_0|=\ell_1+\ell_2$. In addition, $2\ell_1+3\ell_2\leq n $. In order to show the result, there are three cases to be discussed.

{\bf Case 1.} $n=3k$.

 From the fact $2\ell_1+3\ell_2\leq 3k$ we arrive at $\gp(W_n)=|S|=\ell_1+2\ell_2=\frac{2}{3}(2\ell_1+3\ell_2)-\frac{1}{3}\ell_1\leq \frac{2}{3}(3k)-\frac{1}{3}\ell_1=2k-\frac{1}{3}\ell_1$.

On the other hand, $R=\{v_2,v_3,v_5,v_6,v_8,v_9,\cdots,v_{3k-1}, v_{3k}\}$ is a general position set of $W_n$ and $|R|=2k$. Hence, $\gp(W_n)\geq 2k$. That is to say, $2k\leq 2k-\frac{1}{3}\ell_1$, which implies $\ell_1=0.$ Consequently, $\ell_2=k$ and $W_n[S]$  is made up of $2k$ $K_2.$

{\bf Case 2.} $n=3k+1$.

 Using the fact $2\ell_1+3\ell_2\leq 3k+1$, we deduce that $\gp(W_n)=|S|=\ell_1+2\ell_2=\frac{2}{3}(2\ell_1+3\ell_2)-\frac{1}{3}\ell_1\leq \frac{2}{3}(3k+1)-\frac{1}{3}\ell_1=2k+\frac{2}{3}-\frac{1}{3}\ell_1$.

 Now let $R=\{v_2,v_3,v_5,v_6,v_8,v_9,\cdots,v_{3k-1}, v_{3k}\}$. As in Case 1,  $R$ is a general position set of $W_n$ and $|R|=2k$. So, $\gp(W_n)\geq 2k$.

 Hence, $2k\leq 2k+\frac{2}{3}-\frac{1}{3}\ell_1<2k+1.$ Namely, $\gp(W_n)=2k.$ Moreover, $R$ is a $\gp$-set of $W_n.$

{\bf Case 3.} $n=3k+2$.

 From the fact $2\ell_1+3\ell_2\leq 3k+2$ we obtain that  $\gp(W_n)=|S|=\ell_1+2\ell_2=\frac{2}{3}(2\ell_1+3\ell_2)-\frac{1}{3}\ell_1\leq \frac{2}{3}(3k+2)-\frac{1}{3}\ell_1=2k+\frac{4}{3}-\frac{1}{3}\ell_1$.

  Take $R=\{v_2,v_3,v_5,v_6,v_8,v_9,\cdots,v_{3k-1}, v_{3k},v_{3k+2}\}$. It is  not difficult to verify that $R$ is a general position with $|R|=2k+1$. So $\gp(W_n)\geq 2k+1.$

  Hence, $2k+1\leq 2k+\frac{4}{3}-\frac{1}{3}\ell_1<2k+2,$ which implies $\gp(W_n)=2k+1.$ Furthermore, $W_n[S]$ consists of one $K_1$ and $2k$ $K_2$. Evidently, $R$ is a $\gp$-set of $W_n.$

Combining with the above three cases, the result follows.
\qed
\end{proof}

{\bf Acknowledgments} We first would particularly like to thank Professor Klav\v{z}ar for his meaningful and helpful suggestions.  In addition, The authors acknowledge the financial support from National Natural Science Foundation of China (Grant
Nos.11401348 and 11561032),and Shandong Provincial Natural Science Foundation(No. ZR201807061145).


\begin{thebibliography}{99}

\bibitem{1}
B.S. Anand, S.V. Ullas Chandran, M, Changat, S. Klav\v{z}ar, E. J. Thomas, Characterizaation of general position sets and its applications to cographs and bipartite graphs, Appl. Math. Comput. 359(2019)84--89.


\bibitem{2}
J.A. Bondy, U.S.R. Murty, Graph Theory, in: GTM, Springer, 2008.


\bibitem{3}
H.E. Dudeney, Amusements in Mathematics, Nelson, Edinburgh. 1917.

\bibitem{4}
A. Flammenlamp, Progress in the no-three-in-line problem. J. Comb. Theory Ser. A 60(1992)305--311.

\bibitem{5}
V. Froese, I.Kanj, A. Nichterlein, R. Niedermeier, Finding points in general position, Int. J. Comput. Geom. $\&$ Appl. 27(04)(2017) 277--296.

\bibitem{15}
I. Gutman, S.C. Li, W. Wei, Cacti with $n$ vertices and $t$ cycles having extremal Wiener index, Discret. Appl. Math. 232(2017)189--200.

\bibitem{6}
R.R. Hall, T.H. Jackson, A. Sudbery, K. Wild, Some advances in the no-three-in-line problem. J. Comb. Theory Ser. A. 18(1975)336--341.



\bibitem{7}
S. Klav\v{z}ar, I.G. Yero, The general position problem and strong resolving graphs, Open Math. 17 (2019) 1126--1135.

\bibitem{18k}
S. Klav\v{z}ar, B. Patk\'{o}s, G. Rus, I. G. Yero, On general position sets in Cartesian grids, arXiv [math.CO] ( 25 Jul 2019).

\bibitem{8}
C.Y. Kn, K. B. Wong, On no-three-in-line problem on $m$-dimensional torus, Graphs Combin. 34(2018)355--364.


\bibitem{9}
P.Manuel, S. Klav\v{z}ar, A general position problem in graph theory, Bull. Aust. Math. Soc. 98(2018)339--350.



\bibitem{10}
A. Misiak, Z. Stepie\'{n}, A. Szymaszkiewicz, L. Szymaszkiewicz, L. Szymaszkiewicz, M. Zwierzchowski, A note on the no-three-in-line problem on a torus. Discrete Math. 339(2016)217--221.

\bibitem{17p}
B. Patk\'{o}s, On the general position problem on Kneser graphs,
arXiv:1903.08056v2 [math.CO] (21 Jul 2019).

\bibitem{11}
M. Payne, D. R. Wood, On the general position subset selection problem, SIAM J. Discrete Math. 27(2013)1727--1733.

\bibitem{12}
A. Por, D. R. Wood, No-Three-in-Line-in-3D, Algorithmica 47(2007)481--488.

\bibitem{16}
H.M.A. Siddiqui, M. Imran, Computing the metric dimension of wheel related graphs, Appl. Math. Comput. 242(2014)624--632.


\bibitem{13}
M. Skotnica, No-three-in-line problem on a torus:periodicity, Discrete Math. 342 (2019) 111611, 13 pp.

\bibitem{14}
S.J. Wang, On extremal cacti with respect to the Szeged index, Appl. Math. Comput. 309(2017)85--92.

\end{thebibliography}
\end {document}